\documentclass[11pt]{amsart}
\newlength{\basicwidth}
\newlength{\shortbasicwidth}
\newlength{\basicheight}
\numberwithin{equation}{section}

\begin{document}
\title{{Identities involving Bernoulli and Euler polynomials}}
\maketitle

\markboth{\rm \centerline{ H. Alzer and S. Yakubovich}}{}

\markright{\rm \centerline{Bernoulli and  Euler polynomials}}

\vspace{0.5cm}
\begin{center}
HORST ALZER$^a$  \, and  \, SEMYON YAKUBOVICH$^b$
\footnote{
The work of this author  was partially supported by CMUP (UID/MAT/
00144 /2013), which is funded by FCT (Portugal) with national (MEC), European structural funds through the programs FEDER under the partnership agreement PT2020, and Project STRIDE-NORTE-01-0145-FEDER- 000033, funded by ERDF-NORTE 2020.}
\end{center}

\vspace{0.8cm}
\begin{center}
$^a$
Morsbacher Stra$\ss$e 10, 51545 Waldbr\"ol, Germany\\
\emph{email:} \tt{h.alzer@gmx.de}
\end{center}

\vspace{0.3cm}
\begin{center}
$^b$
Department of  Mathematics,
Faculty of Sciences,
University of Porto,
Campo Alegre st. 687,
4169-007 Porto,
Portugal\\
\emph{e-mail:} \tt{syakubov@fc.up.pt}
\end{center}

\vspace{2.8cm}
{\bf{Abstract.}} We present various identities involving the classical Bernoulli and Euler polynomials. Among others, we prove that 
$$
\sum_{k=0}^{[n/4]}(-1)^k {n\choose 4k}\frac{B_{n-4k}(z) }{2^{6k}}
=\frac{1}{2^{n+1}}\sum_{k=0}^{n} (-1)^k
\frac{1+i^k}{(1+i)^k}
{n\choose k}{B_{n-k}(2z)}
$$
and
$$
\sum_{k=1}^{n} 2^{2k-1} {2n\choose 2k-1} B_{2k-1}(z)
=
\sum_{k=1}^n k \,  2^{2k} {2n\choose 2k} E_{2k-1}(z).
$$
Applications of our results lead to formulas for Bernoulli and Euler numbers, like, 
for instance,
$$
n E_{n-1} =\sum_{k=1}^{[n/2]} \frac{2^{2k}-1}{k} (2^{2k}-2^n){n\choose 2k-1} B_{2k}B_{n-2k}.
$$

\vspace{1cm}

\vspace{0.2cm}
\vspace{4mm}

\section {Introduction}

\noindent
The classical Bernoulli and Euler polynomials $B_n(z)$ 
and  $E_n(z)$ $(n=0,1,2,...; z\in\mathbb{C})$ are defined by
\begin{equation}
\frac{t e^{zt}}{e^t-1}=\sum_{n=0}^\infty B_n(z)\frac{t^n}{n!}
\quad
(|t|<2\pi)
\quad\mbox{and}
\quad
\frac{2 e^{zt}}{e^t+1}=\sum_{n=0}^\infty E_n(z)\frac{t^n}{n!} \quad 
(|t|<\pi),
\end{equation}
respectively.
Both functions play an important role in various branches of mathematics, like, for example,
number theory, analysis, combinatorics, statistics and, moreover,  
they  have interesting applications in physics. The first few polynomials are
$$
B_0(z)= 1, \quad  B_1(z)= z- {1\over 2}, \quad  B_2(z)=  z^2- z+ {1\over 6}, 
\quad  B_3(z)=  z^3- {3\over 2} z^2+ \frac{1}{2}z,
$$
$$
E_0(z)= 1, \quad  E_1(z)= z- {1\over 2}, \quad  E_2(z)=  z^2- z, 
\quad  E_3(z)=  z^3- {3\over 2} z^2+ \frac{1}{4}.
$$

The numbers $B_n=B_n(0)$ and $E_n=2^n E_n(1/2)$ are called Bernoulli and Euler numbers.
In particular, we find the values   $B_0=1,\ B_1= -1/2, B_2= 1/6, B_4= -1/30$ and $B_n=0$ for  odd  $n \ge 3$. Furthermore, for all $n\geq 1$,  
 $$
 B_{2n}=(-1)^{n-1}\frac{(2n)!}{2^{2n-1}\pi^{2n}}\zeta(2n),
 $$
 where $\zeta$ denotes the Riemann zeta function.
 All Bernoulli numbers are rational and all Euler numbers are integers. 
 We have $E_0=1$, $E_2=-1$, $E_4=5$ and $E_n=0$ for  odd $n\geq 1$. 
 The Euler numbers can be expressed in terms of Bernoulli polynomials. Indeed, for $n\geq 1$ we have
 $$
 E_{2n}=-\frac{4^{2n+1}}{2n+1} B_{2n+1}(1/4).
 $$
 
 Explicit formulae for $B_n(z)$ and $E_n(z)$ are given by
 $$
 B_n(z)=\sum_{k=0}^n {n\choose k} B_k z^{n-k}=
 \sum_{k=0}^n\frac{1}{k+1}\sum_{j=0}^k(-1)^j {k\choose j}(z+j)^n
 $$
 and
$$
 E_n(z)=\sum_{k=0}^n {n\choose k} \frac{E_k}{2^k}\Bigl( z-\frac{1}{2}\Bigr)^{n-k}=
 \sum_{k=0}^n\frac{1}{2^k}\sum_{j=0}^k(-1)^j {k\choose j}(z+j)^n.
 $$
 The Bernoulli and Euler polynomials and numbers have been studied thoroughly by many authors. 
 Collections of their main  properties
can be found in the books  Apostol \cite[chapter 12]{A}, Erd\'elyi et al. \cite[chapter 1]{erd}, Milovanovi\'c et al. \cite[chapter 1]{MMR}. We also refer to a bibliography of Bernoulli numbers published by Dilcher and Slavutskii \cite{DS}.
Among the huge number of papers on this subject we mention just a few recently published articles:
 Agoh \cite{Ag}, Agoh  and  Dilcher \cite{AD}, Luo \cite{L}, Merca \cite{M},
 Yakubovich \cite{Y}.

 Here,  we list some properties which will be employed below.
 The 
 identities
 \begin{equation}
 \frac{n+1}{2^{n+1}}E_n(2z)=B_{n+1}(z+1/2)-B_{n+1}(z)
\end{equation}
\and
\begin{equation}
\frac{n}{2} E_{n-1} (z)=B_n(z)-2^n B_n(z/2)
\end{equation}
provide elegant connections between the Bernoulli and Euler polynomials.
Moreover,
we have the addition formulae
\begin{equation}
B_n(x+y)=   \sum_{k=0}^{n} {n \choose k} B_k(x) y^{n-k},
\quad
E_n(x+y)=   \sum_{k=0}^{n} {n \choose k} E_k(x) y^{n-k}
\end{equation}
and the multiplication formula
\begin{equation}
B_n(mx)=m^{n-1} \sum_{k=0}^{m-1} B_n\Bigl(x+\frac{k}{m}\Bigr).
\end{equation}

There exist numerous remarkable  identities involving the Bernoulli and Euler polynomials. 
It is the aim of this paper to continue the work on this subject and to
 present various identities for both functions which we could not locate in the literature. In the next section, we 
deduce six identities involving $B_n(z)$ and $E_n(z)$.  In Section 3, we apply some of our results  to obtain additional identities for these polynomials. Finally, in Section 4, we use the identities given in the previous two sections to deduce formulae for the numbers $B_n$ and $E_n$.

\vspace{0.3cm}
\section{Identities}

In this section, we prove six identities for Bernoulli and Euler polynomials. They have in common that they can be deduced  from
certain elementary identities involving the exponential function and the hyperbolic sine and cosine functions.

\vspace{0.3cm}
\noindent
{\bf{Theorem 2.1.}} \emph{For all nonnegative integers $n$ and complex numbers $z$ we have}
\begin{equation}
\sum_{k=0}^{[n/2]}    4^{n-2k} 
 {n\choose 2k} 
{B_{n-2k}(z)}
=\sum_{j=0}^{n} (-1)^{n-j} 2^j {n\choose j}  {B_j(2z)}=2^n B_n(2z-1/2).
\end{equation}

\vspace{0.3cm}
\begin{proof}
A simple calculation reveals that the following identity holds:
\begin{equation}
\cosh(t/4) \, \frac{t e^{zt}}{e^t-1}=e^{-t/4}\, \frac{(t/2) e^{2z \cdot t/2}}{e^{t/2}-1}.
\end{equation}
Next, we develop both sides of (2.2) in a Taylor series. Let
$$
a_{2k}=\frac{1}{4^{2k} (2k)!}
\quad\mbox{and}
\quad{a_{2k+1}=0} \quad{(k=0,1,2,...)}.
$$
We apply (1.1), the Taylor series for $\cosh$ and $\exp$, and the Cauchy product. Then,
\begin{eqnarray}
\cosh(t/4) \, \frac{t e^{zt}}{e^t-1} & = &
\sum_{n=0}^\infty a_n t^n \, \sum_{n=0}^\infty B_n(z)\frac{t^n}{n!}\\\nonumber
& = & \sum_{n=0}^\infty \sum_{j=0}^n a_j \frac{B_{n-j}(z)}{(n-j)!} t^n\\\nonumber
& = & \sum_{n=0}^\infty \sum_{k=0}^{[n/2]} a_{2k} \frac{B_{n-{2k}}(z)}{(n-2k)!} t^n\\\nonumber
& = & \sum_{n=0}^\infty \sum_{k=0}^{[n/2]}  \frac{B_{n-{2k}}(z)}{4^{2k} (2k)!(n-2k)!} t^n
\end{eqnarray}
and
\begin{eqnarray}
e^{-t/4}\, \frac{(t/2) e^{2z \cdot t/2}}{e^{t/2}-1} & = &
\sum_{n=0}^\infty \frac{(-t/4)^n}{n!}\, \sum_{n=0}^\infty B_n(2z)\frac{(t/2)^n}{n!}\\
& = &
\sum_{n=0}^\infty  \sum_{j=0}^n (-1)^{n-j}
\frac{B_j(2z)}{2^{2n-j} j! (n-j)!}t^n.\nonumber
\end{eqnarray}
From (2.2), (2.3), and (2.4) we  obtain the first identity in (2.1).
The second one follows from (1.4) with $x=2z$ and $y=-1/2$.
\end{proof}

\vspace{0.3cm}
\noindent
{\bf{Theorem 2.2.}} \emph{For all nonnegative integers $n$ and complex numbers $z$ we have}
\begin{equation}
\sum_{k=0}^{[n/2]}  \frac{4^{n-2k}}{2k+1} {n\choose 2k} B_{n-2k}(z)
=\sum_{j=0}^{n} (-1)^{n-j} 2^j 
 {n\choose j} 
{E_j(2z)}=2^n E_n(2z-1/2).
\end{equation}

\vspace{0.3cm}
\begin{proof}
The proof of (2.5) is based upon the identity
\begin{equation}
\frac{\sinh(t/2)}{t/2} \, \frac{2t e^{z\cdot 2t}}{e^{2t}-1}=e^{-t/2}\, \frac{2 e^{2z\cdot t} }{e^{t}+1}.
\end{equation}
 Let
$$
b_{2k}=\frac{1}{2^{2k} (2k+1)!}
\quad\mbox{and}
\quad{b_{2k+1}=0} \quad{(k=0,1,2,...)}.
$$
Using (1.1), the Taylor series for $\sinh$ and $\exp$, and the Cauchy product yields
\begin{eqnarray}
\frac{\sinh(t/2)}{t/2} \, \frac{2t e^{z\cdot 2t}}{e^{2t}-1}
& = &
\sum_{n=0}^\infty b_n t^n \, \sum_{n=0}^\infty B_n(z)\frac{(2t)^n}{n!}\\
& = & \sum_{n=0}^\infty \sum_{j=0}^n b_j \frac{2^{n-j} B_{n-j}(z)}{(n-j)!} t^n
\nonumber \\
& = &
\sum_{n=0}^\infty \sum_{k=0}^{[n/2]} b_{2k} \frac{2^{n-2k}B_{n-{2k}}(z)}{(n-2k)!} t^n \nonumber \\
& = &
\sum_{n=0}^\infty \sum_{k=0}^{[n/2]}  \frac{B_{n-{2k}}(z)}{2^{4k-n} (2k+1)!(n-2k)!} t^n\nonumber
\end{eqnarray}
and
\begin{eqnarray}
e^{-t/2}\, \frac{2 e^{2z \cdot t}}{e^{t}+1} & = &
\sum_{n=0}^\infty \frac{(-t/2)^n}{n!}\, \sum_{n=0}^\infty E_n(2z)\frac{t^n}{n!} \\
& = &
\sum_{n=0}^\infty  \sum_{j=0}^n (-1)^{n-j}
\frac{E_j(2z)}{2^{n-j} j! (n-j)!}t^n.\nonumber
\end{eqnarray}
Applying (2.6), (2.7), and (2.8) leads to the left-hand side  of
 (2.5). 
\end{proof}

\vspace{0.3cm}
\noindent
{\bf{Theorem 2.3.}} \emph{For all nonnegative integers $n$ and complex numbers $z$ we have}
\begin{equation}
2^n B_n(z+1/4)=
\sum_{k=0}^{n} 2^{2k-n} {n\choose k}B_k(z) 
=\sum_{k=0}^{n} {n\choose k}{E_{n-k}(1/2)}{B_k(2z)}
=\sum_{k=0}^{n} {n\choose k}{B_{n-k}(1/2)}{E_k(2z)}.
\end{equation}

\vspace{0.3cm}
\begin{proof}
We have 
\begin{equation}
e^{t/2} \, \frac{2t e^{z\cdot 2t}}{e^ {2t}-1}=\frac{1}{\cosh(t/2)}\frac{t e^{2z\cdot t}}{e^t-1}
=\frac{t/2}{\sinh(t/2)}\frac{2 e^{2z\cdot t}}{e^t +1}.
\end{equation}
Next, we apply (1.1), the Taylor series for $\exp$, the representations
$$
{\frac{1}{\cosh(t)}=\sum_{n=0}^\infty E_n \frac{t^n}{n!},
}
\quad{
\frac{t}{\sinh(t)}=\sum_{n=0}^\infty B_n(1/2) \frac{(2t)^n}{n!}
},
$$
 and the Cauchy product. Then we obtain
$$
e^{t/2} \, \frac{2t e^{z\cdot 2t}}
{e^{2t}-1}=\sum_{n=0}^\infty
\frac{1}{n!}
\sum_{k=0}^n {n\choose k} 2^{2k-n} B_k(z) t^n,
$$
$$
\frac{1}{\cosh(t/2)}
\frac{t e^{2z \cdot t}}{e^t-1}=\sum_{n=0}^\infty
\frac{1}{n!}
\sum_{k=0}^n {n\choose k} E_{n-k}(1/2) B_k(2z) t^n,
$$
$$
\frac{t/2}{\sinh(t/2)}
 \frac{2 e^{2z\cdot t}}{e^t+1}=\sum_{n=0}^\infty
 \frac{1}{n!}
 \sum_{k=0}^n {n\choose k} B_{n-k}(1/2) E_k(2z) t^n.
$$
Using (2.10) gives the second and the third identity in (2.9).
\end{proof}

\vspace{0.3cm}
\noindent
{\bf{Theorem 2.4.}} \emph{For all nonnegative integers $n$ and complex numbers $z$ we have}
\begin{equation}
B_n(z)+ 2
\sum_{k=1}^{[n/2]} {n\choose 2k}\frac{B_{n-2k}(z) }{3^{2k+1}}
=\frac{1}{3^n}\sum_{k=0}^{n}(-1)^{n-k} {n\choose k}{B_k(3z)}=\frac{1}{3^n} B_n(3z-1).
\end{equation}

\vspace{0.3cm}
\begin{proof}
We make use of the identity
\begin{equation}
\frac{1+e^{-t}+e^t}{3} \,
\frac{3t e^{z \, 3t}}{e^{3t}-1}=e^{-t}\, \frac{t e^{3z\cdot t}}{e^t-1}.
\end{equation}
Applying
$$
\frac{1+e^{-t}+e^t}{3} =\frac{1}{3}(1+2\cosh(t))=1+\frac{2}{3}\sum_{k=1}^\infty \frac{t^{2k}}{(2k)!}
$$
and (1.1) yields that the product on the left-hand side of (2.12) is equal to
\begin{equation}
\sum_{n=0}^\infty a_n t^n\,\sum_{n=0}^\infty B_n(z)\frac{3^n}{n!} t^n
=\sum_{n=0}^\infty \sum_{k=0}^n a_k\frac{B_{n-k}(z) 3^{n-k}}{(n-k)!} t^n,
\end{equation}
where
\begin{equation}
a_0=1, \, a_{2k-1}=0,
\,\,  a_{2k}=\frac{2}{3\cdot (2k)!} \,\,  (k\geq 1).
\end{equation}
The right-hand side of (2.12) is equal to
\begin{equation}
\sum_{n=0}^\infty (-1)^n \frac{t^n}{n!}
\,\sum_{n=0,}^\infty B_n(3z)\frac{t^n}{n!} 
=\sum_{n=0}^\infty \sum_{k=0}^n (-1)^{n-k}\frac{B_{k}(3z) }{k!(n-k)!} t^n.
\end{equation}
From (2.13) and (2.15) we obtain for $n\geq 0$:
\begin{equation}
\sum_{k=0}^n a_k\frac{B_{n-k}(z)3^{n-k}}{(n-k)!}
=\sum_{k=0}^n (-1)^{n-k}\frac{B_{k}(3z)}{k!(n-k)!}.
\end{equation}
Using (2.14) we get
\begin{equation}
\sum_{k=0}^n a_k \frac{B_{n-k}(z)3^{n-k}}{(n-k)!} 
=\frac{3^n}{n! } B_n(z)+\frac{2\cdot 3^n}{n!}
\sum_{k=1}^{[n/2]} {n\choose 2k}\frac{B_{n-2k}(z)}{3^{2k+1}}.
\end{equation}
From (2.16) and (2.17) we conclude that the first identity in (2.11) holds.
\end{proof}

\vspace{0.3cm}
\noindent
{\bf{Theorem 2.5.}} \emph{For all nonnegative integers $n$ and complex numbers $a\neq 0$, $z$ we have}
\begin{equation}
a^n
\sum_{k=0}^{n} {n\choose k}\frac{B_{n-k}(z) }{k+1}
=E_n(az)+\frac{1}{2}\sum_{k=1}^{n} {n\choose k}{E_{n-k}(az)}=\frac{1}{2}\bigl( E_n(az)+E_n(az+1)\bigr).
\end{equation}

\vspace{0.3cm}
\begin{proof}
We have
\begin{equation}
\frac{e^{at}-1}{at}\,\frac{ate^{z\, at}}{e^{at}-1}=\frac{e^t+1}{2}\frac{2e^{az\, t}}{e^t+1}.
\end{equation}
Using
$$
\frac{e^x-1}{x}=\sum_{k=0}^\infty \frac{x^k}{(k+1)!}
$$
and (1.1) reveals that the left-hand side of (2.19) is equal to
\begin{equation}
\sum_{k=0}^\infty \frac{a^k t^k}{(k+1)!}\,\sum_{n=0}^\infty B_n(z)\frac{a^n t^n}{n!}=\sum_{n=0}^\infty
\sum_{k=0}^n a^n \frac{B_{n-k}(z)}{(k+1)!(n-k)!}t^n.
\end{equation}
Applying
$$
\frac{e^t+1}{2}=1+\frac{1}{2}\sum_{k=1}^\infty \frac{t^k}{k!}
$$
and (1.1) gives that the right-hand side of (2.19) is equal to
\begin{equation}
\sum_{k=0}^\infty c_k t^k
\,
\sum_{n=0}^\infty E_n(az)\frac{t^n}{n!}
=\sum_{n=0}^\infty \sum_{k=0}^n c_k \frac{E_{n-k}(az)}{(n-k)!}t^n,
\end{equation}
where
$$
c_0=1, \,\, c_k=\frac{1}{2\cdot k!}
\,\, (k\geq 1).
$$
From (2.20) and (2.21) we obtain
\begin{eqnarray}
\frac{a^n}{n!}
\sum_{k=0}^n {n\choose k}\frac{B_{n-k}(z)}{k+1}
& = & \sum_{k=0}^n c_k \frac{E_{n-k}(az)}{(n-k)!}\nonumber \\
& = &
\frac{E_n(az)}{n!}+\frac{1}{2\cdot n!}
\sum_{k=1}^n {n\choose k} E_{n-k}(az).\nonumber
\end{eqnarray}
This leads to the left-hand side in (2.18). The second identity follows from (1.4).
\end{proof}

\vspace{0.3cm}
\noindent
{\bf{Theorem 2.6.}} \emph{For all nonnegative integers $n$ and complex numbers $z$ we have}
\begin{equation}
\sum_{k=0}^{[n/4]}(-1)^k {n\choose 4k}\frac{B_{n-4k}(z) }{2^{6k}}
=\frac{1}{2^{n+1}}\sum_{k=0}^{n} (-1)^k
\frac{1+i^k}{(1+i)^k}
{n\choose k}{B_{n-k}(2z)}.
\end{equation}

\vspace{0.3cm}
\begin{proof}
Let
$$
\phi(z)=\cosh\Bigl(
\frac{1+i}{2}z\Bigl)
\cosh\Bigl(\frac{1-i}{2}z
\Bigr)=\frac{1}{2}(\cosh(z)+\cosh(iz))
=\sum_{k=0}^\infty \frac{z^{4k}}{(4k)!}
$$
and
$$
\sigma(z)=\frac{1}{2}(e^{-z}+e^{-iz})=
\sum_{k=0}^\infty (-1)^k \frac{1+i^k}{2}
\frac{z^k}{k!}.
$$
Then,
\begin{equation}
\phi(w)\, \frac{2(1+i)we^{2(1+i)w\, z}}{e^{2(1+i)w}-1}
=\sigma(w)\,\frac{(1+i)w e^{(1+i)w\, 2z}}{e^{(1+i)w}-1}.
\end{equation}
The expression on the left-hand side of (2.23) is equal to
\begin{equation}
\sum_{k=0}^\infty
\frac{w^{4k}}{(4k)!}\,
\sum_{n=0}^\infty \frac{2^n (1+i)^n}{n!}
B_n(z) w^n
=\sum_{n=0}^\infty L_n(z) w^n,
\end{equation}
where
$$
L_n(z)=
 \sum_{k=0}^n a_k\frac{(2(1+i))^{n-k}}{(n-k)!}B_{n-k}(z)
 $$
with
$$
a_{4k}=\frac{1}{(4k)!}, \,\,
a_{4k+1}=a_{4k+2}=a_{4k+3}=0
\quad{(k\geq 0)}.
$$
The expression on the right-hand side of (2.23) is equal to
\begin{equation}
\sum_{k=0}^\infty (-1)^k\frac{1+i^k}{2}
\frac{w^{k}}{k!}\,
\sum_{n=0}^\infty \frac{ (1+i)^n}{n!}
B_n(2z) w^n
=\sum_{n=0}^\infty R_n(z) w^n,
\end{equation}
where
$$
R_n(z)=\frac{1}{2}
 \sum_{k=0}^n (-1)^k \frac{(1+i^k)(1+i)^{n-k}}{k!(n-k)!}B_{n-k}(2z).
 $$
 From (2.24) and (2.25) we obtain for $n\geq 0$:
 $$
 L_n(z)=R_n(z).
 $$
 Since
 $$
L_n(z)=\frac{1}{n!} \sum_{k=0}^{[n/4]} 
{n\choose 4k}(2(1+i))^{n-4k} B_{n-4k}(z),
 $$
 we conclude that (2.22) is valid.
\end{proof}

\vspace{0.3cm}
\section{Applications}

We now apply four of the six theorems presented in Section 2 to obtain new identities involving $B_n(x)$ and $E_n(x)$. 
Theorem 3.3 (below) offers two identities only involving  Bernoulli polynomials, whereas the identities given in the other three theorems
reveal connections between Bernoulli and Euler polynomials. We begin with an application of Theorem 2.1 which leads to the following result.

\vspace{0.3cm}
{\bf{Theorem 3.1.}} \emph{For all nonnegative integers $n$ and complex numbers $z$ we have}
\begin{equation}
\sum_{k=1}^{n} 2^{2k-1} {2n\choose 2k-1} B_{2k-1}(z)
=
\sum_{k=1}^n k \, 2^{2k} {2n\choose 2k} E_{2k-1}(z)
\end{equation}
\emph{and}
\begin{equation}
\sum_{k=0}^{n} 2^{2k} {2n+1\choose 2k} B_{2k}(z)
=
\sum_{k=0}^n  (2k+1) 2^{2k} 
{2n+1  \choose 2k+1} 
E_{2k}(z).
\end{equation}

\vspace{0.3cm}
\begin{proof}
We set in (2.1) $z=x+iy$ ($x,y\in\mathbb{R}$), $n=2N$ and use (1.4). This leads to $L=R$ with
$$
L=
\sum_{k=0}^{N} (-1)^{N-k} 4^{-2k} {2N \choose 2k}\  \sum_{j=0}^{2(N-k)}  {2(N-k) \choose j} i^{-j} B_{j}\left(x\right) y^{ 2(N-k) -j}
$$
and
$$
R=
 \sum_{j=0}^{2N}   (-1)^j  2^{j-4N} {2N\choose j}  \sum_{m=0}^{j}  { j \choose m} B_{m}\left(2 x\right) (2yi)^{j-m}.
$$
Next,
 we split the interior sum of $L$ and the exterior and interior  sums of $R$  into sums with  even and odd indices. This gives
$$
L=
\sum_{k=0}^{N}  (-1)^{N-k} 4^{-2k} {2N\choose 2k}\Bigl[  \sum_{m=0}^{N-k}  (-1)^m {2(N-k) \choose 2m}  B_{2m}(x)  y^{ 2(N-k -m)}
$$
$$
+  i \sum_{m=1}^{N-k}  (-1)^m {2(N-k) \choose 2m-1}  B_{2m-1}(x)
 y^{ 2(N-k -m)+1} \Bigr]
$$
and
$$
R=
\sum_{k=0}^{N}  2^{2k-4N} {2N\choose 2k}\Bigl[  \sum_{j=0}^{k}  (-1)^{k-j}{2k \choose 2j}  B_{2j}(2x)  (2y)^{ 2(k-j)}
$$
$$
+  i \sum_{j=1}^{k}  (-1)^{k-j} {2k \choose 2j-1}  B_{2j-1}(2x)
 (2y)^{ 2(k-j)+1} \Bigr]
$$
$$
+
\sum_{k=1}^{N}   2^{2k-4N-1} {2N\choose 2k-1}\Bigl[ - \sum_{j=1}^{k}  
(-1)^{k-j}  {2k-1 \choose 2j-1}  B_{2j-1}(2x)  (2y)^{ 2(k-j)}
$$
$$
+  i \sum_{j=0}^{k-1}  (-1)^{k-j} {2k -1 \choose 2j}  B_{2j}(2x)
 (2y)^{ 2(k-j)-1} \Bigr].
$$
We now consider  the real  parts of $L$ and $R$. 
$$
\Re{L}=
\sum_{k=0}^{N} \sum_{m=0}^{N-k}  (-1)^{N-k-m}  4^{-2k} {2N\choose 2k} {2(N-k) \choose 2m}  B_{2m}\left(x\right) y^{ 2(N-k -m)}
$$
and
$$
\Re{R}=
   \sum_{k=0}^{N} \sum_{j=0}^{k}     (-1)^{k-j} 2^{-4N+4k-2j} {2N\choose 2k}  {2 k \choose 2j} B_{2j} (2x) y^{ 2(k-j)}
$$
$$
-   \sum_{k=1}^{N} \sum_{j=1}^{k}   (-1)^{k-j}  2^{-4N+4k-2j-1} {2N\choose 2k-1}    {2 k-1 \choose 2j-1} B_{2j-1}\left(2 x\right) y^{ 2(k-j)}.
$$
Using the summation formulae
$$
\sum_{k=0}^N\sum_{m=0}^{N-k} a(k,m)=
\sum_{k=0}^N\sum_{m=0}^{k} a(N-k,k-m),
$$
$$
\sum_{k=0}^N\sum_{j=0}^{k} a(k,j)=
\sum_{k=0}^N\sum_{m=0}^{k} a(k,k-m),
\quad
\sum_{k=1}^N\sum_{j=1}^{k} a(k,j)=
\sum_{k=0}^{N-1}\sum_{m=0}^{k} a(k+1,k+1-m)
$$
and the binomial identities
\begin{equation}
{2N\choose 2k}
{2k\choose 2m}
=
{2N\choose 2m}
{2(N-m)\choose 2(N-k)},
\quad
{2N\choose 2k+1}
{2k+1\choose 2m}
=
{2N\choose 2m}
{2(N-m)\choose 2(N-k)-1}
\end{equation}
yields
\begin{equation}
4^{2N} \Re{L}=
\sum_{k=0}^{N} \sum_{m=0}^{k}  (-1)^{m}  4^{2k} {2N\choose 2m} {2(N-m) \choose 2(N-k)}  B_{2(k-m)}\left(x\right) y^{2m}
\end{equation}
and
\begin{equation}
4^{2N} \Re{R}=
\sum_{k=0}^{N} \sum_{m=0}^{k}  (-1)^{m}  4^{k+m} {2N\choose 2m} {2(N-m) \choose 2(N-k)}  B_{2(k-m)}\left(2x\right) y^{2m}
\end{equation}
$$
+2
\sum_{k=0}^{N-1} \sum_{m=0}^{k}  (-1)^{m+1}  4^{k+m} {2N\choose 2m} {2(N-k) \choose 2(N-k)-1}  B_{2(k-m)+1}\left(2x\right) y^{2m}.
$$
From (1.5) we obtain
\begin{equation}
B_{2(k-m)}(2x)=2^{2(k-m)-1}\bigl[ B_{2(k-m)}(x)+B_{2(k-m)}(x+1/2)\bigr].
\end{equation}
Since $\Re{L}=\Re{R}$, we conclude from (3.4), (3.5) and (3.6) that $U=V$, where
$$
U=
\sum_{k=0}^{N} \sum_{m=0}^{k}  (-1)^{m}  4^{2k-1} {2N\choose 2m} {2(N-m) \choose 2(N-k)} \bigl[B_{2(k-m)}(x)- B_{2(k-m)}(x+1/2)\bigr] y^{2m}
$$
and 
$$
V=
\sum_{k=0}^{N-1} \sum_{m=0}^{k}  (-1)^{m+1}  4^{k+m} {2N\choose 2m} {2(N-m) \choose 2(N-k)-1}  B_{2(k-m)+1}(2x) y^{2m}.
$$
Applying
\begin{equation}
\sum_{k=0}^n\sum_{m=0}^k a(k,m)=\sum_{k=0}^n \sum_{m=0}^{n-k} a(k+m,k)
\end{equation}
and (1.2) gives
$$
U=
\sum_{k=0}^{N} (-1)^{k+1} 2^{4k-1} {2N\choose 2k}
\sum_{m=1}^{N-k}  m  2^{2m} {2(N-k)\choose 2m}  E_{2m-1}(2x) y^{2k}
$$
and
$$
V=
\sum_{k=0}^{N-1} (-1)^{k+1}2^{4k-1}  {2N\choose 2k}\sum_{m=0}^{N-k-1}  2^{2m+1} {2(N-k)\choose 2m+1}   B_{2m+1}(2x) y^{2k}.
$$
Next, we compare the coefficients of $y^{2k}$ $(k=0,1,...,N)$ and find 
\begin{equation}
\sum_{m=1}^{N-k}  m 2^{2m}  {2(N-k)\choose 2m}   E_{2m-1}(2x)=
\sum_{m=0}^{N-k-1}  2^{2m+1} {2(N-k)\choose 2m+1}  B_{2m+1}(2x).
\end{equation}
Finally,
we
replace $N-k$ by $n$,  $2x$ by $z$ and simplify. Then, (3.8) leads to (3.1).

If we consider the identity $\Im{L}=\Im{R}$ and use similar arguments as above then we obtain (3.2). We omit the details.
\end{proof}

\vspace{0.3cm}
Next, we  apply  Theorems 2.2 and 2.4 to prove 
 four identities which are given in the following two theorems.

\vspace{0.3cm}
{\bf{Theorem 3.2.}} \emph{For all nonnegative integers $n$ and complex numbers $z$ we have}
\begin{equation}
\sum_{k=0}^{n} \frac{4^{2k}}{2(n-k)+1} {2n\choose 2k} B_{2k}(z)
=
\sum_{k=0}^n   2^{2k} {2n\choose 2k} E_{2k}(2z)
-\sum_{k=0}^{n-1}   2^{2k+1} {2n\choose 2k+1} E_{2k+1}(2z)
\end{equation}
\emph{and}
\begin{equation}
\sum_{k=0}^{n} \frac{4^{2k+1}}{2(n-k)+1} {2n+1\choose 2k+1} B_{2k+1}(z)
=
\sum_{k=0}^n   2^{2k+1} {2n+1\choose 2k+1} E_{2k+1}(2z)
-\sum_{k=0}^{n}   2^{2k} {2n+1\choose 2k} E_{2k}(2z).
\end{equation}

\vspace{0.3cm}
\begin{proof}
We replace  in (2.5) $z$ by  $x+iy$ $(x,y\in\mathbb{R})$ and $n$ by  $2N$. Then,  we study the real parts of the left-hand side and the right-hand side, respectively. An application of  (1.4) gives
$$
\Re{L}=\sum_{k=0}^N  \frac{4^{2k}}{2(N-k)+1   } {2N\choose 2k} 
\sum_{m=0}^k (-1)^{k-m} {2k\choose 2m} B_{2m}(x) y^{2(k-m)}
$$
and
$$
\Re{R} = 
\sum_{k=0}^N  { 2^{2k} } {2N\choose 2k} 
\sum_{m=0}^k (-1)^{k-m} 2^{2(k-m)}  {2k\choose 2m} E_{2m}(2x) y^{2(k-m)}
$$
$$
-
\sum_{k=0}^{ N-1}{ 2^{2k+1} } {2N\choose 2k+1} 
\sum_{m=0}^k (-1)^{k-m} 2^{2(k-m)} {2k+1\choose 2m+1} E_{2m+1}(2x)
 y^{2(k-m)}.
$$
We use the summation formula
\begin{equation}
\sum_{k=0}^N\sum_{m=0}^k a(k,m)=\sum_{k=0}^N\sum_{m=k}^N a(m,m-k)
\end{equation}
and obtain
$$
\Re{L} =
\sum_{k=0}^N \sum_{m=k}^N \frac{(-1)^{k}}{2(N-m)+1} 4^{2m}{2N\choose 2m}
{2m\choose 2k} B_{2(m-k)}(x) y^{2k},
$$
$$
\Re{R} =
\sum_{k=0}^N \sum_{m=k}^N (-1)^{k} 2^{2m+2k}{2N\choose 2m}
{2m\choose 2k} E_{2(m-k)}(2x) y^{2k}
$$
$$
-\sum_{k=0}^{N-1} \sum_{m=k}^{N-1} (-1)^{k} 2^{2m+2k+1}{2N\choose 2m+1}
{2m+1\choose 2k} E_{2(m-k)+1}(2x) y^{2k}.
$$
Now, we  compare the coefficients of $y^{2k}$ $(k=0,...,N)$, apply 
(3.3), and simplify. This yields
$$
\sum_{m=k}^N \frac{4^{2m}}{2(N-m)+1}{2(N-k)\choose 2(N-m)} B_{2(m-k)}(x)
$$
$$
=\sum_{m=k}^N 2^{2m+2k} {2(N-k)\choose 2(N-m)} E_{2(m-k)}(2x)
-
\sum_{m=k}^{N-1} 2^{2m+2k+1}{2(N-k)\choose 2(N-m)-1} E_{2(m-k)+1}(2x).
$$
We use
\begin{equation}
\sum_{m=k}^N a(m)=\sum_{m=0}^{N-k} a(m+k).
\end{equation}
Then, we set $N-k=n$  and simplify. This yields (3.9).

If we study the imaginary parts of both sides of (2.5) with $z=x+iy$, $n=2N$ and use the method of proof given above, then we obtain (3.10).
\end{proof}

\vspace{0.3cm}
{\bf{Theorem 3.3.}} \emph{For all nonnegative integers $n$ and complex numbers $z$ we have}
\begin{equation}
\sum_{k=0}^{n}  {2n\choose 2k} B_{2k}(3z)
-
\sum_{k=0}^{n-1}    {2n\choose 2k+1} B_{2k+1}(3z)
= 2
\sum_{k=0}^{n-1}   3^{2k-1} {2n\choose 2k} B_{2k}(z)+3^{2n} B_{2n}(z)
\end{equation}
\emph{and}
\begin{equation}
\sum_{k=0}^{n}  {2n+1\choose 2k+1} B_{2k+1}(3z)
-
\sum_{k=0}^{n}    {2n+1\choose 2k} B_{2k}(3z)
= 2
\sum_{k=0}^{n-1}   3^{2k} {2n+1\choose 2k+1} B_{2k+1}(z)+3^{2n+1} B_{2n+1}(z).
\end{equation}

\vspace{0.3cm}
\begin{proof}
We set in (2.11)  $z=x+iy$ $(x,y\in \mathbb{R})$,  $n=2N$
and make use of (1.4). Then we obtain for the real part of the left-hand side,
$$
\Re{L}=\sum_{k=0}^N {2N\choose 2k} (-1)^{N-k} B_{2k}(x) y^{2(N-k)}
$$
$$
+2\sum_{k=1}^N {2N\choose 2k} \frac{1}{3^{2k+1}}\sum_{m=0}^{N-k} {2(N-k)\choose 2m} (-1)^{N-k-m} B_{2m}(x) y^{2(N-k-m)},
$$
and for the real part of the right-hand side we obtain
$$
\Re{R}=
\frac{1}{9^N}\Bigl[
\sum_{k=0}^N {2N\choose 2k} \sum_{m=0}^{k} 3^{2(k-m)}(-1)^{k-m}
{2k\choose 2m}  B_{2m}(3x) y^{2(k-m)}
$$
$$
-\sum_{k=0}^{N-1}  {2N\choose 2k+1} \sum_{m=0}^{k} 3^{2(k-m)}(-1)^{k-m}
{2k+1\choose 2m+1}  B_{2m+1}(3x) y^{2(k-m)}
\Bigr].
$$
Using the summation formulae (3.11) and
$$
\sum_{k=1}^N\sum_{m=0}^{N-k} a(k,m)=\sum_{k=0}^{N-1}\sum_{m=k}^{N-1}a(N-m,m-k)
$$
 we find
$$
\Re{L}=\sum_{k=0}^N {2N\choose 2k} (-1)^{k} B_{2(N-k)}(x) y^{2k}
+2\sum_{k=0}^{N-1}\sum_{m=k}^{N-1}  {2N\choose 2m}  {2m\choose 2k} (-1)^{k} 3^{2(m-N)-1} B_{2(m-k)}(x) y^{2k}
$$
and 
$$
\Re{R}=
\frac{1}{9^N}\Bigl[
\sum_{k=0}^N \sum_{m=k}^N {2N\choose 2m} {2m\choose 2k}
 3^{2k}(-1)^{k}
 B_{2(m-k)}(3x) y^{2k}
$$
$$
-\sum_{k=0}^{N-1} \sum_{m=k}^{N-1}
 {2N\choose 2m+1} {2m+1\choose 2k}
  3^{2k}(-1)^{k}
  B_{2(m-k)+1}(3x) y^{2k}
\Bigr].
$$
Next, we compare the coefficients of $y^{2k}$ $(k=0,1,...,N)$, make use of (3.3), and simplify. This yields
$$
 B_{2(N-k)}(x)+2\sum_{m=k}^{N-1} {2(N-k)\choose 2(N-m)}
3^{2(m-N)-1} B_{2(m-k)}(x)
$$
$$
=\frac{1}{9^N}
\Bigl[
\sum_{m=k}^{N}  {2(N-k)\choose 2(N-m)}
3^{2k} B_{2(m-k)}(3x)
-
\sum_{m=k}^{N-1} {2(N-k)\choose 2(N-m)-1}
3^{2k} B_{2(m-k)+1}(3x)
\Bigr].
$$

Finally, we apply (3.12),
set $N-k=n$ and simplify again. This gives (3.13).

 If we consider the impaginary parts of both sides of (2.11) with $z=x+iy$, $n=2N$,  then we obtain (3.14).
\end{proof}

\vspace{0.3cm}
We conclude this section with an application of Theorem 2.6 which leads us to  four closely related identities for certain alternating sums involving Bernoulli and Euler polynomials.

\vspace{0.3cm}
{\bf{Theorem 3.4.}}
\emph{For all positive integers $n$ and complex numbers $z$ we have}
\begin{equation}
\sum_{k=1}^n (-1)^{k} 2^{2k+1} k {4n\choose 4k} E_{4k-1}(z)
\end{equation}
$$
=
\sum_{k=1}^n (-1)^{k} 2^{2k-1} {4n\choose 4k-1} B_{4k-1}(z)
-
\sum_{k=1}^n (-1)^{k} 2^{2k-2}  {4n\choose 4k-3} B_{4k-3}(z),
$$
\begin{equation}
\sum_{k=1}^n (-1)^{k} 2^{2k+1} (2k-1) {4n-2\choose 4k-2} E_{4k-3}(z)
\end{equation}
$$
=
\sum_{k=1}^{n-1} (-1)^{k} 2^{2k+1} {4n-2\choose 4k-1} B_{4k-1}(z)
+
\sum_{k=1}^n (-1)^{k} 2^{2k}  {4n-2\choose 4k-3} B_{4k-3}(z),
$$
\begin{equation}
\sum_{k=0}^n (-1)^{k} 2^{2k} (4k+1) {4n+1\choose 4k+1} E_{4k}(z)
\end{equation}
$$
=
\sum_{k=0}^{n-1} (-1)^{k} 2^{2k+1} {4n+1\choose 4k+2} B_{4k+2}(z)
+
\sum_{k=0}^n (-1)^{k} 2^{2k}  {4n+1\choose 4k} B_{4k}(z),
$$
\begin{equation}
\sum_{k=0}^n (-1)^{k} 2^{2k+1} (4k+3) {4n+3\choose 4k+3} E_{4k+2}(z)
\end{equation}
$$
=\sum_{k=0}^{n} (-1)^{k} 2^{2k+1} {4n+3\choose 4k+2} B_{4k+2}(z)
-
\sum_{k=0}^{n} (-1)^{k} 2^{2k} {4n+3\choose 4k} B_{4k}(z).
$$

\vspace{0.3cm}
\begin{proof}
We replace in (2.22) $n$ by $4N$. Then we obtain for the expression on the left-hand side:
$$
L=\frac{(-1)^N}{2^{6N}}
\sum_{k=0}^N (-1)^k 2^{6k} {4N\choose 4k} B_{4k}(z),
$$
and for the expression on the right-hand side we get
$$
R=\frac{(-1)^N}{2^{6N}}
\Bigl[
\sum_{k=0}^N (-1)^k 4^{k} {4N\choose 4k} B_{4k}(2z)
$$
$$
-\frac{1}{2}
\sum_{k=1}^N (-1)^k 4^{k} {4N\choose 4k-1} B_{4k-1}(2z)
-
\sum_{k=0}^{N-1} (-1)^k 4^{k} {4N\choose 4k+1} B_{4k+1}(2z)\Bigr].
$$
We set $z=x+iy$ $(x,y\in \mathbb{R})$ and apply (1.4). Then we get for the real parts of $L$ and $R$:
$$
(-1)^N 2^{6N}
\Re{L}=
\sum_{k=0}^N\sum_{m=0}^{2k} (-1)^{k+m}4^{3k}
{4N\choose 2m}{4N-2m\choose 4k-2m} B_{2m}(x) y^{4k-2m},
$$
$$
(-1)^N 2^{6N}
\Re{R}=
\sum_{k=0}^N \sum_{m=0}^{2k} (-1)^{k+m}4^{3k-m}
{4N\choose 2m}{4N-2m\choose 4k-2m} B_{2m}(2x) y^{4k-2m}
$$
$$
+\frac{1}{2}
\sum_{k=0}^{N-1}  \sum_{m=0}^{2k+2} (-1)^{k+m}4^{3k+3-m}
{4N\choose 2m-1}{4N-2m+1\choose 4k-2m+4} B_{2m-1}(2x) y^{4k+4-2m}
$$
$$
-\sum_{k=0}^{N-1}  \sum_{m=0}^{2k} (-1)^{k+m}4^{3k-m}
{4N\choose 2m+1}{4N-2m-1\choose 4k-2m} B_{2m+1}(2x) y^{4k-2m}.
$$
Applying (1.3) we conclude that $\Re{L}=\Re{R}$ is equivalent to $U=V$, where
$$
U=
\sum_{k=0}^N\sum_{m=0}^{2k} (-1)^{k+m}4^{3k-m}m
{4N\choose 2m}{4N-2m\choose 4k-2m} E_{2m-1}(2x) y^{4k-2m},
$$
$$
V=\sum_{k=0}^{N-1} \sum_{m=0}^{2k} (-1)^{k+m}4^{3k-m}
{4N\choose 2m+1}{4N-2m-1\choose 4k-2m} B_{2m+1}(2x) y^{4k-2m}
$$
$$
-\frac{1}{2}
\sum_{k=0}^{N-1} \sum_{m=0}^{2k+2} (-1)^{k+m}4^{3k+3-m}
{4N\choose 2m-1}{4N-2m+1\choose 4k-2m+4} B_{2m-1}(2x) y^{4k+4-2m}.
$$
We split the interior sums into two sums and find
$$
U=
\sum_{k=0}^N\sum_{m=0}^{k} (-1)^{k}4^{k+2m} 2(k-m)
{4N\choose 4m}{4(N-m)\choose 4(k-m)} E_{4(k-m)-1}(2x) y^{4m}
$$
$$
+\sum_{k=0}^{N-1}\sum_{m=0}^{k} (-1)^{k}4^{k+2m+2}(2(k-m)+1)
{4N\choose 4m+2}{4(N-m)-2\choose 4(k-m)+2} E_{4(k-m)+1}(2x) y^{4m+2}
$$
and $V=V_1+V_2$ with
$$
V_1
=\sum_{k=0}^{N-1}  \sum_{m=0}^{k} (-1)^{k}4^{k+2m}
{4N\choose 4m}{4(N-m)\choose 4(k-m)+1} B_{4(k-m)+1}(2x) y^{4m}
$$
$$
+\frac{1}{2}
\sum_{k=0}^{N}  \sum_{m=0}^{k} (-1)^{k}4^{k+2m}
{4N\choose 4m}{4(N-m)\choose 4(k-m)-1} B_{4(k-m)-1}(2x) y^{4m},
$$
and
$$
V_2
=\sum_{k=0}^{N-2}  \sum_{m=0}^{k} (-1)^{k}4^{k+2m+2}
{4N\choose 4m+2}{4(N-m)-2\choose 4(k-m)+3} B_{4(k-m)+3}(2x) y^{4m+2}
$$
$$
+\frac{1}{2}
\sum_{k=0}^{N-1}  \sum_{m=0}^{k} (-1)^{k}4^{k+2m+2}
{4N\choose 4m+2}{4(N-m)-2\choose 4(k-m)+1} B_{4(k-m)+1}(2x) y^{4m+2}.
$$
Next, we use the summation formula 
(3.7) and compare the coefficients of $y^{4k}$ and $y^{4k+2}$. This leads to the identities
$$
\sum_{m=1}^{N-k} (-1)^m 2^{2m+1} m {4(N-k)\choose 4m} E_{4m-1}(2x)
$$
$$
=\sum_{m=1}^{N-k} (-1)^m 2^{2m-1} {4(N-k)\choose 4m-1} B_{4m-1}(2x)
+\sum_{m=0}^{N-k} (-1)^m 2^{2m} {4(N-k)\choose 4m+1} B_{4m+1}(2x)
$$
and
$$
\sum_{m=0}^{N-k-1} (-1)^m 2^{2m} (2m+1) {4(N-k)-2\choose 4m+2} 
E_{4m+1}(2x)
$$
$$
=\sum_{m=0}^{N-k-2} (-1)^m 2^{2m}  {4(N-k)-2\choose 4m+3} B_{4m+3}(2x)
+\sum_{m=0}^{N-k-1} (-1)^m 2^{2m-1} {4(N-k)-2\choose 4m+1} B_{4m+1}(2x).
$$
Finally, we replace $N-k$ by $n$, $2x$ by $z$, and simplify. This gives (3.15) and (3.16).

 If we study the imaginary parts of $L$ and $R$ and use the same method of proof as above, then we get (3.17) and (3.18).
\end{proof}

\vspace{0.3cm}
\section{Bernoulli and Euler numbers}

In this section, we use some of our theorems to find identities involving Bernoulli and Euler numbers. We only present a few illustrative examples. Many more identities can be obtained easily. To prove our results we apply several known formulae which are valid for all natural numbers $n$:
$$
B_n(1)=(-1)^n B_n,
\quad
B_{2n}(1/6)=\frac{1}{2} (1-2^{1-2n})(1-3^{1-2n})B_{2n},
\quad
B_{2n}(1/3)=-\frac{1}{2}(1-3^{1-2n})B_{2n},
$$
$$
B_n(1/2)=(2^{1-n}-1)B_n,
\quad
B_n(1/4)=2^{-n}(2^{1-n}-1)B_n-n 4^{-n} E_{n-1},
$$
$$
E_{n-1}(0)=\frac{2}{n}(1-2^n) B_n,
\quad
E_{2n}(1/6)=2^{-2n-1}(1+3^{-2n})E_{2n}.
$$
We get the following two formulae by setting  $z=1/2$ and $z=1/4$, respectively, in (2.5).
$$
E_n=\sum_{k=0}^{[n/2]} \frac{2^{n-2k}}{2k+1} (2-2^{n-2k}){n\choose 2k} B_{n-2k},
$$
$$
\sum_{k=1}^n \frac{2^{2k}-2}{2(n-k)+1}{2n\choose 2k} B_{2k}=\frac{1}{2n+1}.
$$
From  (2.9) with $z=0$ and $z=1/4$ we obtain
$$
\sum_{k=0}^n 2^k (2^k-E_{n-k}) {n \choose k} B_{k}=0,
$$
$$
(2n+1) E_{2n} =\sum_{k=0}^{n} \frac{4^{k+1}-1}{k+1} (2^{2k+2}-2^{2n+1}){2n+1\choose 2k+1} B_{2k+2}B_{2n-2k},
$$
$$
\sum_{k=0}^n (2^k-2){n\choose k} B_k (E_{n-k}-1)=\sum_{k=1}^n k{n\choose k} E_{k-1}.
$$

Using (2.22) with $z=0$ gives for odd integers $n\geq 3$:
\begin{displaymath}
\sum_{k=0}^n (-1)^k \frac{1+i^k}{(1+i)^k} {n\choose k} B_{n-k}
= \left \{ \begin{array}{ll}
{(-1)^{(n+3)/4}  }{2^{(3-n)/2}} n , & \textrm{if $n \equiv 1 \, (\mbox{mod} \, 4)$,}\\
0, & \textrm{if $n \equiv 3 \, (\mbox{mod} \, 4)$.} 
\end{array} \right.
\end{displaymath}

Next, we set in (3.1) and (3.2) $z=0$. Then we get the closely related formulae
$$
\sum_{k=1}^n 2^{2k}(2^{2k}-1) {2n\choose 2k} B_{2k}= 2n
\quad\mbox{and}
\quad
\sum_{k=0}^n 2^{2k}{2n+1\choose 2k} B_{2k} =2n+1.
$$
The identity
$$
\sum_{k=0}^n (2^{2k-1}-1)(1-3^{1-2k}){2n+1\choose 2k} B_{2k}
=\sum_{k=0}^n (k+1/2) (1+3^{-2k}) {2n+1\choose 2k+1} E_{2k}
$$
follows from (3.2) by setting $z=1/6$. An application of
(3.10) with $z=0$ yields
$$
\sum_{k=1}^n \frac{2^{2k}}{k} (2^{2k}-1){2n-1\choose 2k-1} B_{2k}=2.
$$
 From (3.13) with  $z=1/3$ we obtain a representation of $B_{2n}$:
$$
\frac{1}{2} (3^{2n}-1) B_{2n}=\sum_{k=0}^{n-1}(1-3^{2k-1}){2n\choose 2k}B_{2k}.
$$
Finally, we set $z=0$ in (3.15), (3.16), (3.17) and (3.18), respectively. This leads to
$$
\sum_{k=0}^{n}(-1)^k 2^{2k}(4^{2k}-1) {4n\choose 4k} B_{4k}=2n,
$$
$$
\sum_{k=0}^{n}(-1)^k 2^{2k+1}(4^{2k+1}-1) {4n+2\choose 4k+2} B_{4k+2}=2n+1,
$$
$$
\sum_{k=0}^{n-1} (-1)^{k} 2^{2k+1} {4n+1\choose 4k+2} B_{4k+2}
+
\sum_{k=0}^{n}(-1)^k 2^{2k} {4n+1\choose 4k} B_{4k}=
4n+1,
$$
$$
\sum_{k=0}^{n} (-1)^{k} 2^{2k+1} {4n+3\choose 4k+2} B_{4k+2}
=
\sum_{k=0}^{n} (-1)^{k} 2^{2k} {4n+3\choose 4k} B_{4k}.
$$

\vspace{5mm}

\end{document}